\documentstyle[12pt]{article}
\normalbaselines
\begin{document}
\vbox {\vspace{6mm}}
\def\WP{W(\q,\phi)}
\def\WPa{W(\q,\phi, a)}
\def\FPa{F_a(\phi)}
\def\FP{F(\phi)}
\def\0{\Pi_0}
\def\APR{A^{\phi}(\R)}
\def\OP{\Omega^{\phi}}
\def\AP{A^{\phi}}
\def\hP{\h^{\phi}}
\def\hPR{\hP(\R)}
\def\HPR{\HP(\R)}
\def\HP{H^{\phi}}

\def\ll{\lambda}
\def\o{\omega}
\def\O{\Omega}
\def\W{{\bf W}(\phi)}
\def\r{{\bf r}}
\def\g{{\bf g}}
\def\Wr{{\bf W}(\r,\phi)}
\def\Wb{{\bf W}_{\Phi}}
\def\d{{\bf d}}
\def\R{{\bf R}}
\def\C{{\bf C}}
\def\gC{{\bf g}({\C})}
\def\1C{{\bf g}_1({\C})}
\def\2C{{\bf g}_2({\C})}
\def\dC{{\bf d}({\C})}
\def\dR{{\bf d}(\R)}
\def\gR{{\bf g}(\R)}
\def\lq{\bar\q}
\def\lqq{^{\lq}}
\def\w{{\bf w}}
\def\h{{\bf h}}
\def\t{{\bf t}}
\def\a{\alpha}
\def\E{E_{\a}}
\def\F{E_{-\a}}
\def\H{H_{\a}}
\def\u{{\bf u}}
\def\q{\sigma}
\def\e{\varepsilon}
\def\s{s^{\e}}
\def\th{\theta}
\def\k{{\bf k}}
\def\m{{\bf m}}
\def\b{\beta}
\def\n{{\bf n}}
\def\N{{\bf N}}
\def\p{{\bf p}}
\def\l{{\bf l}}
\def\GC{G(\C)}
\def\GR{G(\R)}
\def\ss{{{\bf g}^*}}
\def\ssC{{\bf g}^*({\C})}
\def\G{\Gamma}
\def\qq{^{\q}}
\def\lh{\overline\h}
\def\lr{\overline\r}
\def\1{\Pi_1}
\def\2{\Pi_2}
\def\WR{W(\R)}
\def\gr{\g_{\R}}
\def\lfR{\l_1(\R)}
\def\lsR{\l_2(\R)}
\def\nfR{\n_1^+(\R)}
\def\nsR{\n_2^-(\R)}
\def\lrfR{\lr_1(\R)}
\def\lrsR{\lr_2(\R)}
\begin{center}
{\bf MANIN TRIPLES OF REAL SIMPLE LIE ALGEBRAS. PART 2} 
\end{center}
\begin{center}
A.PANOV \\
Samara State University\\
ul.Akad.Panlova 1,
Samara, 443011,Russia\\
panov@info.ssu.samara.ru
\end{center}
{\bf Abstract.} 
We complete the study of Manin triples of real simple Lie algebras.
In the Part 2 of the article we classify the Manin triples 
$(\gR,W,\gR\oplus\gR)$ (case 2 of the doubles) 
up to weak and gauge equivalence.

\begin{center}
{\bf  0.Introduction} 
\end{center}

First we recall the main definitions of the Part 1.\\
{\bf Definition 1}. Let $\g_1,\g_2,\d$ be  Lie algebras over a field $K$
and let $Q$ be a symmetric  nondegenerate bilinear form on $\d$.
A triple $(\g_1, \g_2,\d)$ is called a Manin triple if 
 $Q(x,y)$ is $ad$-invariant and $\d$ is a direct sum of 
maximum isotropic subspaces $\g_1, \g_2$ . \\
{\bf Definition 2}. We say that two Manin triples 
$(\g,W,\d)$ and $(\g,W',\d)$
are weak equivalent if there exists  an element $a$ in the adjoint group
$D$ of the double $\d$ such that $ W'= Ad_a(W)$.\\
{\bf Definition 3}. We say that two Manin triples 
$(\g,W,\d)$ and $(\g,W',\d)$
are gauge equivalent if there exists  an element $a$ in the adjoint group
$D$ of double  $\d$ such that $ W'=Ad_a(W)$ and $Ad_a (\g) =\g$.

We study Manin triples as follows.
First for given $\g$ we find all the doubles $\d$ and describe 
all forms $Q(\cdot,\cdot)$ on $\d$ such that $\g$ is an isotropic subspace
with respect to this form. Then fixing $\d$ and $Q$ we study 
complementary subalgebras $W$ such that $(\g,W,\d)$ forms a Manin triple. 

We are going the classify Manin triples of real simple Lie algebras up to weak
and gauge equivalence. It is known that the double $\d$ of complex simple Lie 
algebra $\g$ is isomorphic to one of following Lie algebras
$\g\otimes A_i(\C)$ , where $A_1(\C) = \C[t]/t^2$,
$A_2(\C) = {\C[t]}/{(t^2-1/4)}$.  
In Part 1 of this paper we  study
the action of conjugation $\q$ of $\C/\R$ on complex doubles and obtain
that the double $\dR$ of real simple Lie  algebra $\gR$ is isomorphic to
one of following algebras $\g\otimes A_i(\R)$, where
$A_1(\R) = \R[t]/t^2$,
$A_2(\R) = \R[t]/(t^2-1/4)$, $A_3(\R) = \C$.

So, according to the double, there exist 3 types of Manin triples.
In Part 1 we consider Manin triples $(\gR,W,\g)$  (i.e.Case 3).
In this paper we study Manin triples $(\g,W,\g\oplus\g)$ (i.e. Case 2.)\\

\begin{center}
{\bf On Manin triples $(\g,W,\g\oplus\g)$}
\end{center}

Let $\g$ be a semisimple 
Lie algebra over the field  $\C$ of complex numbers.
Let $G$ be the adjoint group of $\g$. 
Denote by $\Delta $ ( resp. $\Delta ^+$)  
the system of roots (resp.of positive roots) and by 
$\Pi=\{\a_1,\ldots ,\a_n\}$ 
the sistem of simple roots. Let $K(\cdot,\cdot)$ be a Killing form. 
Consider the Weyl basis : $$H_k=H_{\a_k}, \a_k\in \Pi; E_{\a}, E_{-a}, 
\a\in\Delta^+$$  
$$ [H,\E]=\a(H)\E , \quad [\E ,\F]=H_{\a} ,\quad $$
$$\a(H)=K(H_{\a}, H),\quad K(\E,\F)=1$$
As usial $\h$ is a Cartan subalgebra and $\n^{\pm}$ is upper and lower
nilpotent subalgebras.

Let $(\g,W,\d)$ be a Manin triple with $\d=\g\oplus\g$.
We recall the classification of Manin triples in this case.
Denote  by $j$ the image of $t$ in $A_2(\C)$.
and  $e=\frac{1}{2}+j$, $f=\frac{1}{2}-j$.
The elements $e$, $f$ are ortogonal idempotents, 
i.e. $ef=fe=0$ and $e^2=f^2=1$.
The double $\d = \g e + \g f$.
The simple calculations yield that 
$$Q(ae+bf, a'e+b'f) = K(a,a')-K(b,b')$$ 
where $K(\cdot,\cdot)$ is the Killing form of $\g$.

Let $\1$, $\2$ be the subsets of $\Pi$ and $\phi:\1\mapsto\2$ is the map
such that:\\
i) $\phi$ is a bijection,\\
ii) $(\phi(\a),\phi(\b))=(\a,\b)$,\\
iii) for every $\a\in\1$ there exists $k$  with the property
$\phi(\a),\ldots,\phi^{k-1}(\a)\in\1$ and $\phi^k(\a)\notin\1$.\\
Notations:\\
$\g_i$ the Lie semisimple subalgebra of $\g$ with the
system of simple roots $\Pi_i$, system of positive roots $\Delta_i^+$
and Cartan subalgebra $\h_i$; $G_i$ is the subgroup with Lie algebra $\g_i$.\\
$\phi:\g_1\mapsto\g_2$ isomorphism of semisimple Lie subalgebras defined
via $\phi:\1\mapsto\2$,\\ 
$\p_i^{\pm}$ is upper and lower parabolic subalgebras,
$ \n_i^{\pm}$ is the nilpotent radical of $\p_i^{\pm}$,\\
$\r_i=\g_i+\h$ is the reductive part of $\p_i^{\pm}$, 
$R_i$ is Lie subgroup with Lie algebra  $\r_i$ and $R' = R_1\bigcap R_2$.\\ 
Let $\lh_1$, $\lh_2$ be the subalgebras of $\h$ such that 
$\lh_i^\perp\subset\lh_i$ and $\h_i\subset\lh_i$.\\
Denote $\l_i =\lh_i^\perp = Ker K(\cdot,\cdot)\vert_{\lh_i}$, 
$\lr_i = \g_i + \lh_i$\\
Suppose that  $\phi$ is extended  to the isomorphism
$$\phi :\lr_1/\l_1\mapsto\lr_2/\l_2\eqno (1.1)$$
{\bf Definition 1.1}. 
The following condition we shall call the condition iv):\\
iv) $\phi$ preserve $K(\cdot,\cdot)$ and has no fixed points in $\h$.\\
( the last means that there is no $x\in\lr_1\bigcap\lr_2$, $x\ne0$ such that 
$\phi(x)=x(mod(\l_2))$.\\
Consider
$$ \W = \l_1 e + \n_1^+ e +\sum_{x\in\lr_1}\C(xe + \phi(x)f)
+ \n_2^ - f +\l_2 f\eqno(1.2)$$
{\bf Theorem 1.2}( [BD],[ST],[CGR],[S]). Let $\g$ be a semisimple Lie algebra.
Let $\phi$ satisfy i), ii), iii), iv). Then \\
1) $(\g, \W, \g\oplus\g))$ is a Manin triple;\\
2) Every Manin triple $(\g, W',\g\oplus\g)$ is gauge 
equivalent to $(\g, \W, \g\oplus\g)$.\\ 

\begin{center}
{\bf On Manin triples $(\gR,W,\gR\oplus\gR)$}
\end{center}

Let $\gR$ be a real simple Lie algebra.
It is well known that $\gR$ is a real form of a complex simple Lie 
algebra $\g$ or coinsides with the realification of $\g$ (see [VO],[GG]).
It was proved in [CGR] and [CH] that  Manin triples 
$(\gR,W,\gR\oplus\gR)$ 
correspond
to the solutions $\tilde Q:\gR\mapsto\gR$ of Modified Yang-Baxter equation
$$[ \tilde Q x, \tilde Q y] - \tilde Q[\tilde Q x,y] - \tilde Q[x, \tilde Q y] 
= \lambda [x,y]$$
with $\lambda <0$.
A solution exists whenever the simple real Lie algebra
$\gR$ has no black roots in Satake-diagram [CGR]. 
It takes place for algebras:\\
L1) $\gR$ is a split real form complex simple Lie algebra ;\\
L2) $\gR = \g_\R$ is the realification of $\g$;\\
L3) $ \gR= su(p,p),su(p,p+1)$;\\
L4) $\gR= so(p,p+2)$;\\
L5) $\gR = EII$;\\
Let $\gR$ be a real simple Lie algebra of L1-L5.\\
Denote by $\G=\{1,\q : \q^2=1\}$ the Galois group $\C/\R$. There exists
a semilinear action of $\G$ on $\gC=\g\otimes\C$ such that $\gR = \gC^\G$.
There exist the subalgebras $\h$, $\n^\pm$ in $\g\C$ such that 
$\qq\n^+= \n^+$ , $\qq \n^-= \n^-$, $\qq\h= \h$ (see [CGR]).
We are fixing further this subalgebras.\\
We have the action of $\q$ on the roots ${\qq\a}(h)= \overline{\a({\qq h})}$, 
$h\in\h$, $\a\in \Delta$.
In all cases  L1-L5 $\q(\Pi)=\Pi$. For the acton of $\q$ on $\Pi$ see
Satake-diagrams [VO].

Let $(\gR, W,\dR)$ be a Manin triple with
$\dR=\gR\oplus\gR$ (i.e. Case 2).
Then $(\gC, W\otimes\C, \dR\otimes\C)$ is Manin triple over $\C$.
Clearly, $\dC = \dR\otimes\C = \gC\oplus\gC$. 
The action of $\q$ on $\dC$ is as follows $\qq e= e$ and $\qq f = f$.

Recall that $\gC=\g$ in the case $\gR$ is a real form of 
a complex simple Lie algebra $\g$ and $\gC=\g+\g$ in the case
$\gR=\g_\R$ is the realification of $\g$.
To simplify the notations, we shall refer $\gC$ as 
$\g$ till the end of this section.

According to Theorem 1.2,  $W\otimes\C$ coinsides with some $\q$-invariant
subalgebra $Ad_g\W$, where $g\in G$.In this section we are going to study  
$\q$-invariant  subalgebras in the orbit $Ad_g\W$ of the adjoint group $G$.
The list of real forms of this subalgebras $ W = (Ad_g\W)^\Gamma$
provides the classification of Manin triples  $(\gR,W,\gR\oplus\gR)$.

Consider $$ ^\q(Ad_g\W) =Ad_g\W\eqno (2.1)$$
It follows
$${^\q(\W)} = Ad_p\W,\eqno (2.2)$$
where $p={^\q g^{-1}} g$. Note that $^\q p = p^{-1}$.  
Replacing $\W$ by (1.2) in (2.1), we get the equality of two  subalgebras
$$\begin{array}{c}^\q\l_1e + ^\q \n_1^+ e +
\sum_{x\in\lr_1}\C(^\q xe + ^\q(\phi(x))f)
+^\q \n_2^ - f + {^\q\l_2}f =   
Ad_p(\l_1)e +\\ 
Ad_p(\n_1^+)e+
\sum_{x\in\lr_1}\C(Ad_p(x)e + Ad_p(\phi(x))f)
+ Ad_p(\n_2^ -) f + Ad_p(\l_2)f\end {array}\quad  \eqno (2.3)$$
Denote $\m_1=\g_1\bigcap\n_2^-$, $\m_2=\g_2\bigcap\n_1^+$.\\
{\bf Lemma 2.1}. We save the notation of section 1. 
Let $\phi$ satisfy i), ii), iii), iv) and (2.3) hold.
Then the following statements are true.\\
1) The conjugation $\q$ saves $\l_1$, $\r_1$, $\lr_1$, $\g_1$, $\m_1$ 
and $\l_2$, $\r_2$, $\lr_2$, $\g_2$, $\m_2$.
This defines the maps 
$\q :\lr_1/\l_1\mapsto\lr_1/\l_1$,
$\q :\lr_2/\l_2\mapsto\lr_2/\l_2$;\\
2) $p\in P_1^+\bigcap P_2^-$; there exists the decomposition 
$p=  m_1qm_2$ with 
$m_1\in\exp(\m_1)$, $ m_2 \in exp(\m_2)$, $q\in R'$.
The elements $m_1$, $m_2$, $q$ obey the following relations
$m_1m_2=m_2m_1$, ${\qq m_1}= q^{-1}m_1^{-1}q$, ${\qq m_2}= qm_2^{-1}q^{-1}$, 
${\qq q}= q^{-1}$;\\
3) Denote $ p_1 = m_1q$ and $p_2 = qm_2$. Then $\qq p_1=p_1^{-1}$ and
$\qq p_2=p_2^{-1}$ ;\\
4) $Ad_{p_1}(\lr_1)=\lr_1$, $Ad_{p_1}(\l_1)= \l_1$,
$Ad_{p_2}(\lr_2) = \lr_2$,
$Ad_{p_2}(\l_2)= \l_2$. This defines
$Ad_{p_1}:\lr_1/\l_1\mapsto\lr_1/\l_1$,
$Ad_{p_2}:\lr_2/\l_2\mapsto\lr_2/\l_2$
 and
$$
\phi(Ad_{p_1}^{-1}{^\q(x))} = Ad_{p_2}^{-1}{^\q\phi(x)}
\eqno (2.4)$$
holds for $x\in\lr_1/\l_1$.\\
{\bf Proof}.
Consider the intesection of left and right sets of (2.3) with $\g e$ 
and $\g f$.\\
We have 
$$\left \{ \begin{array}{rcl}^\q(\l_1+\n_1^+)&=& Ad_p(\l_1+\n_1^+)\\
^\q(\l_2+\n_2^-)& = &Ad_p(\l_2+\n_2^-)\end{array}\right.\eqno (2.5)$$
Compairing the radicals of this subalgebras, we have

$$\left \{ \begin{array}{rcl}^\q\n_1^+&=& Ad_p\n_1^+\\
^\q\n_2^-& = &Ad_p\n_2^-\end{array}\right.\eqno (2.6)$$
Since ${\qq \n_1^+} = \n_1^+$ and ${\qq \n_2^-} = \n_2^-$, then 
$$\left \{ \begin{array}{rcl}\n_1^+&=& Ad_p\n_1^+\\
\n_2^-& = &Ad_p\n_2^-\end{array}\right.\eqno (2.7)$$
This implies $$p\in P_1^+\bigcap P_2^-.$$ 
The element $p$ can be uniquely decomposed 
$$p = m_1qm_2,\eqno (2.8)$$
where $m_1\in exp(\m_1)$, $m_2\in exp(\m_2)$, $q\in R'$.\\
First we are going to prove statement 2) of the Lemma.
Choose the elements $\mu_1\in\m_1$ and $\mu_2\in\m_2$ such that
$m_1= exp(\mu_1)$, $m_2= exp(\mu_2)$. Up to definition,
$$
\mu_1=\sum_{\a\in\Delta_1-\Delta_2}\xi_\a E_{-\a},\qquad
\mu_2=\sum_{\b\in\Delta_2-\Delta_1}\eta_\b E_{\b}
$$
Suppose that $[E_{-\a}, E_\b]= cE_\gamma$ and $c\ne 0$. Then $-\a+\b=\gamma$ 
is a root. If $\gamma >0$, then $\b=\a+\gamma$ and there exists the simple roots
$\1-\2$ in decomposition of $\b$ in the sum of simple roots. This contradicts
to $\b\in\Delta_2$. If $\gamma <0$ , then $\a=\b-\gamma$ and similarly 
we get a contradiction. Hence $[E_{-\a}, E_\b]=0$. 
Therefore, $[\mu_1,\mu_2]=0$ and $m_1m_2= m_2m_1$.

Since ${^\q p}= p^{-1}$ , then
$${\qq m_1}{\qq q} {\qq m_2}= m_2^{-1}q^{-1}m_1^{-1}= 
m_2^{-1}(q^{-1}m_1^{-1}q)q^{-1}=$$ 
$$ (q^{-1}m_1^{-1}q)m_2^{-1}q^{-1}= 
 (q^{-1}m_1^{-1}q)q^{-1}(qm_2^{-1}q^{-1}).$$
Thus, 
${\qq m_1}= q^{-1}m_1^{-1}q$, ${\qq m_2}= qm_2^{-1}q^{-1}$, 
${\qq q}= q^{-1}$. This proves 2). One can deduce 3) from 
2) by direct calculations.

Since $\g_2\bigcap\l_2 = 0$, then $Ad_{p_2}\l_2 = \l_2$.
Therefore,
$$Ad_p(\l_2+\n_2^-) = Ad_p\l_2 + Ad_p\n_2^- = 
Ad_{m_1}\l_2 + \n_2^- =\l_2 +\n_2^-$$ 
Simillarly, $Ad_{p_1}\l_1 = \l_1$ and 
$$Ad_p(\l_1+\n_1^+) = Ad_p\l_1 + Ad_p\n_1^+ = 
Ad_{m_2}\l_1 + \n_1^+ =\l_1 +\n_1^+$$ 
Substituting this to (2.5), we have
$$\left \{ \begin{array}{rcl}^\q(\l_1+\n_1^+)&=& \l_1+\n_1^+\\
^\q(\l_2+\n_2^-)& = &\l_2+\n_2^+\end{array}\right.\eqno (2.9)$$
Since $\qq\h=\h$, then
$$\left \{ \begin{array}{rcl}^\q(\l_1)&=& \l_1\\
^\q(\l_2)& = &\l_2\end{array}\right.\eqno (2.10)$$
The (2.3) can be rewritten as
$$\begin{array}{c}\l_1e + \n_1^+ e +
\sum_{x\in\lr_1}\C(^\q xe + ^\q(\phi(x))f)
+\n_2^ - f + \l_2f =\\   
\l_1e + \n_1^+e + \sum_{x\in\lr_1}\C(Ad_p(x)e + Ad_p(\phi(x))f)
+ \n_2^- f + \l_2f\end {array}\eqno (2.11)$$
We stress that\\
 $Ad_px= Ad_{p_1}x(mod(\n_1^+))$, 
$Ad_{p_1}x\in\lr_1\subset\r_1$
and $Ad_p\phi(x) =Ad_{p_2}\phi(x)(mod(\n_2^+))$, $Ad_{p_2}\phi(x) \in\lr_2
\subset\r_2$.
We obtain
$$ 
\sum_{x\in\lr_1}\C(^\q xe + ^\q(\phi(x))f)
\subset \l_1e + \n_1^+e +
\sum_{x\in\lr_1}\C(Ad_{p_1}(x)e + Ad_{p_2}(\phi(x))f)
+ \n_2^ - f + \l_2f \eqno (2.12)$$
Since $p_1=h_1g_1$ with $h_1\in H=exp(\h)$, $g_1\in G_1=exp(\g_1)$,
then 
$$\left \{ \begin{array}{rcl}Ad_{p_1}\lr_1&=& \lr_1\\
Ad_{p_2}\lr_2& = &\lr_2\end{array}\right.\eqno (2.13)$$
$$\left \{ \begin{array}{rcl}Ad_{p_1}\l_1&=& \l_1\\
Ad_{p_2}\l_2& = &\l_2\end{array}\right.\eqno (2.14)$$
The (2.12) implies that for every $x\in\lr_1$ there exists
$y\in\lr_1$ such that
$$\left \{ \begin{array}{rcl}^\q(x) &=&Ad_{p_1}y(mod(\l_1))\\
^\q(\phi(x))& = &Ad_{p_2}\phi(y)(mod(\l_2))\end{array}\right.\eqno (2.15)$$
It follows
$$\left \{ \begin{array}{rcl}^\q(\lr_1)&=& \lr_1\\
^\q(\lr_2)& = &\lr_2\end{array}\right.\eqno (2.16)$$
The (2.10) and (2.16) prove the statement 1). 

We have seen above that $Ad_{p_1}$ saves $\lr_1$, $\l_1$ and
$Ad_{p_2}$ saves $\lr_2$, $\l_2$. Finally the (2.15) yields 
$\phi(Ad_{p_1}^{-1}{^\q(x))} = Ad_{p_2}^{-1}{^\q\phi(x)}$
for $x\in\lr_1/\l_1$.  $\Box$\\
{\bf Definition 2.2}. 
The following condition we shall call the condition v*) for
$\phi :\lr_1/\l_1\mapsto\lr_2/\l_2$:\\
v*) $\qq\lr_i=\lr_i$, $\qq \l_i=\l_i$ for $i=1,2$ and 
 $${\qq(\phi(x))} = \phi({\qq x})\eqno(2.17)$$ for $x\in\lr_1/\l_1.$\\
Since $\q:\Pi\mapsto\Pi$, then $\qq\Pi_1=\Pi_1$, $\qq\2=\2$ and 
$\qq\phi(\a)= \phi({\qq\a})$ for $\a\in\1$.\\ 
Notations:\\
$\0=\1\bigcup\2$;\\
$\hP=\{v\in\h\vert \a(v) = \phi(\a){v}, \a\in\1\} = 
\{v\in\h\vert [v,\phi(x)] = \phi([v,x]), x\in\lr_1\}$;\\
The last formula can be rewritten as $ad_v\phi(x)=\phi ad_v(x), x\in\lr_1$
We shall call $\hP$ the stabilizer subalgebra of $\phi$.\\
$\HP=\{p\in H\vert Ad_p\phi(x) = \phi Ad_p(x)\}$. We shall call $\HP$ the 
stabilizer of $\phi$. The Lie algebra of $\HP$ is $\hP$.\\
{\bf Definition 2.3}[BD].Let $\phi:\1\mapsto\2$ satisfy i),ii),iii). 
Let $\a,\b\in\1$. We shall say that $\a < \b$ if there exists integer
$k$ such that $\phi^k(\a)=\b$.\\
Note that $\a < \phi(\a)<\cdots<\phi^k(\a)=\b$.\\
{\bf Definition 2.4}. The set $\a_1<\a_2<\cdots<\a_k$ is called a chain.
We shall denote by $C(\a)$ the maximal chain with $\a$.\\
{\bf Lemma 2.5}. 1) The maximal chains either coinside or don't intersect.;
2) If $\phi$ and $\q$ obey v*) (see Definition 2.2), then the image $\q(C)$ 
of a maximal chain is also a maximal chain.\\
{\bf Proof}.\\
1) Let $C_1= \{\a_1 < \a_2 <\cdots<\a_k\}$ and 
$C_2= \{\b_1 < \b_2 <\cdots<\b_k\}$ are maximal chains. If
$\phi^n(\a) = \phi^m(\b)$ , $n\ge m$ then $\phi^{n-m}(\a)=\b$ and
, therefore, $\b\in C_1$. Hence, $C_2\subset C_1$. Since $C_2$ is maximal,
 $C_1 =C_2$.\\
2) Let $\{\a_1 < \a_2 <\cdots < \a_k\}$ be a maximal chain with 
$\a_{i+1} = \phi(\a_i)$. Denote $\b_i=\q(\a_i)$. We see
$\phi(\b_i) = \phi(\q\a_i)) = \q(\phi(\a_i))) = \q(\a_{i+1}) =\b_{i+1}$.
Hence, $\b_1 < \b_2 <\cdots < \b_k$.$\Box$\\
{\bf Corollary 2.6} The $\0$ is decomposed in maximal chains. The automorphism 
$\q$ acts on the set of maximal chains.\\ 
{\bf Theorem 2.7}. Let $\gR$ be a simple Lie algebra of L1-L5 types.
Suppose that the map $\phi $ satisfies i)-iv)
and  the algebra  $Ad_g\W $ is  $\q$-invariant. Denote as above
$p={\qq g^{-1}}g$ . Then $\phi$ satisfies v*) and  $p\in H^\phi$.\\
{\bf Proof}.
The proof is divided into 2 steps.\\
1) We recall that (2.1) implies (2.2) with $p={^\q g^{-1}}g$. 
We are going to show in this step that
if $Ad_g\W$ is $\q$-invariant then $\phi$ satisfies v*).

It was proved that (2.1) implies (2.4). 
The equality (2.4) can  be rewritten in the form
$$Ad_{p_2}^{-1}{\qq\phi}\qq(Ad_{p_1}x) = \phi(x)\eqno(2.18)$$
for $x\in\lr_1/\l_1$.\\
Decompose $\1$ and $\2$ in connected components:
$$\1 = \bigcup_{i=1}^k \Pi_{1i}$$
$$\2 = \bigcup_{i=1}^k \Pi_{2i}$$
Then $\g_1$ and $\g_2$ are decomposed in direct sums:
$$\g_1 = \sum_{i=1}^k \g_{1i}$$
$$\g_2 = \sum_{i=1}^k \g_{2i}$$
We may assume that $\phi(\g_{1i})=\g_{2i}$.

As we saw in Lemma 2.1 $^\q\g_1=\g_1$, $^\q\g_2=\g_2$.
Since $p_1\in R_1$, then 
$Ad_{p_1}(\g_1) = \g_1$ and 
$Ad_{p_1}$ save the decomposition of $\g_1$ in the some of simple algebras.
Simillarly, $p_2\in R_2$, then 
$Ad_{p_2}(\g_2) = \g_2$  and 
$Ad_{p_2}$ save the decomposition of $\g_2$ in the some of simple algebras.
Denote $\psi={\qq\phi}\qq$.
According to (2.18), the maps $\psi$ and $\phi$ gave the same permutation of
simple components . Thus 
$\psi(\g_{1i})=\g_{2i}$. Therefore, $\phi^{-1}\psi(\g_{1i})=\g_{1i}$
and $\phi^{-1}\psi(\Pi_{1i})= \Pi_{1i}$. 
 
Denote $\phi^{-1}\psi=s_{1i}$.
Then $\psi=\phi s_{1i}$. Substituting $\phi s_{1i}$
for $\psi$ in (2.18), we get  
$$Ad_{p_2}^{-1}\phi s_{1i}(Ad_{p_1}x) = \phi(x)\eqno(2.19),$$

$$\phi^{-1}Ad_{p_2}^{-1}\phi s_{1i}(Ad_{p_1}x) = x\eqno(2.19)$$
for $x\in g_{1i}$.
Denote $v_2\in\r_2$ such that $p_2=exp(v_2)$. 
Then
$$ \phi^{-1}Ad_{p_2}^{-1}\phi(x) = \phi^{-1}exp (-ad_{v_2})\phi(x)= 
exp(- ad _{\phi^{-1}v_2})x$$
Therefore, $\phi^{-1}Ad_{p_2}^{-1}\phi$ is an inner automorphism of
$\g_{1i}$. We denote it by $Ad_q$. The equality
$Ad_qs_{1i}Ad_{p_1}=1$ implies that $s_{1i}$ is inner automorphism of
$\g_{1i}$. Since $s_{1i}= \phi^{-1}\psi$ save $\Pi_{1i}$, then $s_{1i}=1$.
Finally, we get $\phi^{-1}\psi=1$ and $\psi=\phi$.
This proves (2.17) for $x\in\g_1$.
If $x$ belongs to orthogonale complement of $\g_1$ in $\lr_1$ , then
(2.4) is equivalent to (2.17).
Consequently, (2.17) holds for all $x\in\lr_1/\l_1$.\\
2) The goal of this step is to prove that $p_1=p_2=p\in H^\phi$.
As we prove in step 1), (2.17) and (2.4) hold.
Then     
$$\phi(Ad_{p_1}^{-1}{^\q(x))} = Ad_{p_2}^{-1}{^\q\phi(x)} = 
  Ad_{p_2}^{-1}{\phi(\qq x)}$$
We get
$$\phi(Ad_{p_1}x) = Ad_{p_2}\phi(x),\eqno(2.19)$$
for all $x\in\lr_1/\l_1$.

There exists $v_i\in\r_i$ such that $p_i= exp(v_i)$.
According Lemma 2.1, we have
$v_1=\mu_1 +\th + t$, 
$v_2= \mu_2+\th+t$ where $\mu_i\in\m_i$, 
$\th\in(\g_1\bigcap\g_2)\bigcap(\n^+\oplus \n^-)$, $t\in \h$.
We shall prove further that $\mu_1=\mu_2=\th=0$.
 
We can decompose $t=t_i+t_i'$ where $t_i\in\h_i=\g_i\bigcap\h$, $t_i'\in\h$
and $[t_i',\g_i]= 0$. Denoting $g_i=exp(\mu+\th+t_i)$, $h_i=exp(t_i')$, we get
$p_i= g_ih_i= h_ig_i$. 

The (2.19) implies 
$\phi(Ad_{g_1}x) = Ad_{g_2}\phi(x)$. 
Let $G_i$ be the subgroup of $G$ with Lie algebra $\g_i$.
We can extend $\phi$ to the isomorphism of
the groups $\phi:G_1\mapsto G_2$. Then 
$Ad_{\phi(g_1)}\phi(x) = Ad_{g_2}\phi(x)$. Hence, $\phi(g_1)=g_2$ .
Therefore, $\phi(\mu_1+\th + t_1)= 
\mu_2+ \th + t_2$. The element $\th$ is a sum $\th= \th^++ \th^-$ where
$\th^\pm\in\g_1\bigcap\g_2\bigcap\n^\pm$. We obtain
$$\phi(\mu_1)+\phi(\th^+) + \phi(\th^-) + \phi(t_1)= 
\mu_2+ \th^++ \th^-+ t_2\eqno(2.20)$$
Recall that the map $\phi$ is defined via the map $\phi:\1\mapsto\2$.
Thus, $\phi(\n^+)\bigcap\g_1\subset\n^+\bigcap\g_2$, 
$\phi(\n^-)\bigcap\g_1\subset\n^-\bigcap\g_2$, $\phi(\h_1)\subset\h_2$.
The equality (2.20) holds if and only if
$$\left \{ \begin{array}{rcl}\phi(\th^+)&=& \mu_2+\th^+\\
 \phi(\mu_1)+\phi(\th^-)& = &\th^-\\
\phi(t_1)&=&t_2\end{array}\right.\eqno (2.21)$$
Suppose that $\th^+\ne0$. Denote $\Delta_i^+$ the system if positive roots
generated by $\Pi_i$. Then
$$
\th^+=\sum_{\a\in\Delta_1^+}\xi_\a E_\a,\qquad
\mu_2= \sum_{\a\in\Delta_2^+-\Delta_1^+}\xi_\b E_\b$$
The equality $\phi(\th^+)=\mu_2+ \th^+$ implies
$$
\sum_{\a\in\Delta_1^+}\xi_\a E_{\phi(\a)} = 
\sum_{\a\in\Delta_2^+-\Delta_1^+}\xi_\b E_\b
+ \sum_{\a\in\Delta_1^+}\xi_\a E_\a\eqno(2.22)
$$
For every simple root $\a\in\1$ we denote $l(\a)$ the length 
of the maximal chain  $C(\a)$ (see definition 2.4) from the 
beginning to $\a$. 
Denote
$$l_0=min\{l(\a)\vert \a\in\Delta_1^+, \xi_\a\ne0\}$$
and $l_0=l(\a_0)$ for some positive root $\a_0\in\Delta_1^+$
The term $\xi_{\a_0}E_{\a_0}$ belongs to the right side of (2.22)
and dose not belong to to left side. A contradiction.

Therefore $\th^+=0$ and $\mu_2=0$.
Similarly we get  $\th^-=0$ and $\mu_1=0$.
This proves $v_1=v_2=t$ , $p_1=p_2=exp(t)=p\in H$
and 
$$ Ad_p\phi = \phi Ad_p\eqno(2.23)$$
This proves the Theorem. $\Box$\\
{\bf Corollary 2.8} Let $\phi:\1\mapsto\2$ satisfy i)-iv).
Consider the follwing conditions:\\
1) $Ad_g\W$ is $\q$-invariant , 
2) $\W$ is $q$-invariant.
Then 1) implies 2) and 2) is equivalent to v*).\\
{\bf Proof}. If $Ad_g\W$ is $\q$-invariant for some $g\in G$, 
then $\phi$ obeies v*). This proves 1)$\mapsto$ v*) and 2)$\mapsto$v*).

According to Lemma 2.1,if $Ad_g\W$ is $\q$-invariant, then
 $\qq\l_i= \l_i$, $\qq\n_1^+= \n_1^+$, 
$\qq\n_2^-= \n_2^-$, $\qq\lr_i= \lr_i$ . 
The condition v*) is true and yields   
$$ \qq(xe+\phi(x)f) = 
{\qq x}e + {\qq \phi(x)}f) =   
{\qq x}e + \phi({\qq x})f $$
for $x\in\lr_1/\l_1$. He have
$$\qq \W = {\qq (\l_1 e + \n_1^+ e +\sum_{x\in\lr_1}\C(xe + \phi(x)f)
+ \n_2^ - f +\l_2 f)} = $$
$$\l_1 e + \n_1^+ e +\sum_{x\in\lr_1}\C(xe + \phi(x)f)
+ \n_2^ - f +\l_2 f= \W.$$
This proves 1)$\mapsto$2) and v*)$\mapsto$2).
$\Box$\\

\begin{center}
{\bf Completion of classification}
\end{center}
In this section we are going to classify all Manin triples
$\gR,W,\gR\oplus\gR$ up to the $Ad$-action of the group
$G(\R)=G(\C)^\Gamma$.
 As above $\gC=\gR\otimes\C$. Recall that $\gC=\g$ in the case $\gR$ 
is a real form of a complex simple Lie algebra $\g$ and 
$\gC=\g+\g$ in the case $\gR=\g_\R$ is the realification of $\g$.
Denote $\Pi^c$ the set of simple roots of $\gC$. So $\Pi^c=\Pi$ in the case
of real forms and $\Pi^c=\Pi\bigcup\Pi'$ in the case $\gR$ is the realification 
of $\g$ ( here $\Pi$,$\Pi'$ are sets of simple roots of $\g$). 
Let we have two subsets $\Pi_1^c$, $\Pi_2^c$ in $\Pi^c$.
Let $\phi:\1^c\mapsto\2^c$ satisfy i)-iiv).
The $\W$ is the subalgebra (1.2), constructed by the map $\phi$ for $\gC$. \\
{\bf Proposition 3.1}. 
Let $\gR$ be a real simple Lie algebra if the type L1-L5.
and let the subalgebra $Ad_g\W$ be $\q$-invariant.
We assert 
$(Ad_g\W)^\Gamma = Ad_r(\W)^\Gamma$ for some $r\in\GR$.\\
{\bf Proof}.
Recall that we choose $\q$-invariant subalgebras $\n^\pm$, 
$\h$ in $\gC$. The conjugation $\q$ saves 
the set of simple roots $\Pi$.
According to Theorem 2.7., 
the element $p={\qq g}^{-1}g$ lies in $H^\phi$, $\W$ is also
$\q$-invariant and $\phi$ satisfies conditions v*).
Then $\phi(\qq\a)= {\qq\phi}(\a)$ for $\a\in\1$.

Denote by $v$ the element in $\h$ such that $p=exp(v)$. 
The element $p$ uniquelly determined by the system of numbers
$p_\a= exp(\a(v))$, $\a\in\Pi$. 

As $p= exp(v)$, then $\qq p = exp({\qq v})$ and 
$$({\qq p})_\a = exp \a({\qq v}) = exp\overline{\qq \a(v)} = 
\overline{p_{\q(\a)}}$$ 
Easily
$$(p^{-1})_\a = exp( \a(-v))= p_\a^{-1}$$
The equality $\qq p=p^{-1}$ implies  
$$\overline{p_{\q(\a)}}= p_\a^{-1}\eqno(3.1)$$ 
for all $\a\in\Pi$.
Since $p\in H^{\phi}$, then $Ad_p\phi(x)= \phi(Ad_p(x))$, $x\in\lr_1/\l$.
This is equivalent to $Ad_p\phi(E_\a)= \phi(Ad_p(E_\a))$, $\a\in\1$ .
We obtain
$$p_\a= p_{\phi(\a)}\eqno(3.2)$$
for $\a\in\1$.

By Lemmas 2.5-2.6, $\Pi_0$ is decomposes in maximal chains and $\q$ acts on
the chains by permutations.
It follows from (3.2) that  the system of elements
$p_\a$ is stable on maximal chains in $\Pi_0$. 
There exists a 
system  $h_\a\in\C^*$ which is 
stable on maximal chains and 
$p_\a = \overline{h_{\q(\a)}}^{-1}h_\a$, $\a\in\Pi$.
Therefore, there exists the element $h\in H^\phi$ such that
$p= {\qq h^{-1}}h$.
By definition, $\qq g^{-1}g= p= {\qq h^{-1}h}$.
Hence $\qq(gh^{-1}) = gh^{-1}$ and $gh^{-1}\in G(\R)$. 
Finally $g=rh$, $r\in G(\R)$, $h\in H^\phi$.
We have
$$ Ad_h\W = Ad_h(\l_1 e + \n_1^+ e) +\sum_{x\in\lr_1}\C(Ad_hxe + Ad_h\phi(x)f)
+ Ad_h(\n_2^ - f +\l_2 f) = $$ 
$$ \l_1 e + \n_1^+ e +\sum_{x\in\lr_1}\C(Ad_hxe + \phi(Ad_hx)f)
+ \n_2^ - f +\l_2 f= \W.$$
It follows $Ad_g\W = Ad_r\W$ and 
$(Ad_g\W)^\Gamma = Ad_r(\W)^\Gamma$.$\Box$

Let $\g(\R)$ be a real simple Lie algebra of types L1-L2.
We choose $\q$-invariant subalgebras $\h$, $\n^\pm$ in $\gC$ as above.
The $\h(\R)$ is a real form of $\h(\C)$. The conjugation $\q$ acts onthe set 
of simple roots as in Satake-diagrams of L1-L5.
Let $\1^c$ and $\2^c$ be $\q$-invariant subsets in the set of 
simple roots of $\gC$. 
Let  the map $\phi:\1\mapsto\2$ be $\q$-invariant and 
satisfy the conditions i), ii), iii).\\
Denote by $\g_i(\R)$, $\h_i(\R)$ the real forms of $\g_i(\C)$, $\h_i(C).$
Choose two real subspases $\h_i(\R)$ such that 
$\lh_i(\R)^\perp\subset\h_i(\R)$ and $\h_i(\R)\subset\lh_i(\R)$.\\
Denote $\lr_i(\R)= \g_i(\R)+\lh_i(\R)$, $\l_i(\R)=\lh_i(\R)^\perp.$
One can extend the map $\phi$ to the isomorphism
$$\phi:\g_1(\R)\mapsto\g_2(\R).$$
Suppose $\phi$ is extended to the isomorphism
$$\phi:\lr_1(\R)/\l_1(\R)\mapsto\lr_2(\R)/\l_2(\R),$$
preserving Killing form and having no fixed points.
We shall refer the above condition as condition iv$\R$) 
Consider the subalgebra
$$ \W_{\R} = \lfR+ \n_1(\R)^+ e +\sum_{x\in\lr_1(\R)} \R(xe + \phi(x)f)
+ \n_2(\R)^ - f +\lsR \eqno(3.3)$$
{\bf Theorem 3.2}. 
Let $\phi$ be $\q$-invariant and satisfy i)-iii) and iv$\R$). 
Let  $\gR$ be a simple Lie algebra of type L1-L5. 
Then\\
1) The triple $(\gR,\W_{\R},\gR\oplus\gR)$
is a Manin triple;\\
2) Every Manin triple $(\gR,V,\gR\oplus\gR)$ is $Ad_{G(\R)}$ equivalent to
the triple $(\gR,\W_{\R},\gR\oplus\gR)$.\\
{\bf Proof}. 
Consider the complexifications  $\g_i(\C)$, $\lr_i(\C)$, $\l_i(\C)$ of  
$\g_i(\R)$, $\lr_i(\R)$, $\l_i(\R)$.
We extend the $\phi$ to $\C$-linear $\q$-invariant isomorphism
of $\g_1(\C)$ onto $\g_2(\C)$ and further to
$\lr_1(\C)/\l_1(\C)$ onto $\lr_2(\C)/\l_2(\C)$.
We get $\phi$, satisfying conditions i)-iv), and v*).
Therefore, $\W$ is $\q$-invariant.
The algebra $\W^\Gamma $ is generated by  
$\lfR$, $\lsR$, $\nfR$, $\nsR$
and
the elements
$$ xe+\phi(x)f + \qq(xe+\phi(x)f) = 
(x+ {\qq x})e + (\phi(x)f+ {\qq \phi(x)}f) =   
(x+ {\qq x})e + \phi(x+ {\qq x})f $$
where $x\in\lr_1(\C)/\l_1(\C).$
We have $\W^\Gamma = \W_\R$ and $(\gR,\W_\R,\gR\oplus\gR)$
is a Manin triple. this proves 1).

Let $(\gR,W,\gR\oplus\gR)$ be a Manin triple;
then $(\gC, W\otimes\C,\gC\oplus\gC)$ is also a Manin triple.
The algebra  $W\otimes\C$ is $\q$-invariant and, according to Theorem 1.2, 
$W\otimes\C = Ad_g\W.$ 
Then $\phi$ satisfies v*) (see Theorem 2.7) and we get the isomorphism iv$\R$) 
of the algebras of invariants.
We obtain
$(W\otimes\C)^\Gamma = (Ad_g\W)^\Gamma= Ad_r(\W)^\Gamma= Ad_r\W_\R.$
This proves the Theorem.$\Box$\\
{\bf Corollary 3.3} Let $\gR$ be the split form of a simple complex 
Lie algebra $\gR$. Then $\lfR=\lsR=0$ and
$$ \W_{\R} = \n_1(\R)^+ e +\sum_{x\in\lr_1(\R)} \R(xe + \phi(x)f)
+ \n_2(\R)^ - f \eqno(3.4)$$
All Manin triples 
$(\gR, W, \gR\oplus\gR)$ has the form $W= Ad_r(\W_\R)$ with $r\in\GR$.\\
{\bf Proof}.
The Killing form $K(\cdot, \cdot)$ is positive definite on $\h(\R)$.
Therefore, $\lfR= \lsR = 0$.$\Box$\\

Further we shall specify the statement of Theorem 3.2 in the case
$\gR$ is a realification of $\g$.\\
{\bf Theorem 3.4} Let $\g$ be a complex simple Lie algebra and $\gR=\g_\R$
is the realification of $\g$.
as in section 1 let $\phi\1\mapsto\2$ satisfy i)-iii).
We extend $\phi$ to $\R$-isomorphism
$\phi:\g_1\mapsto\g_2$ such that $\phi$ is $\C$ linear or $\C$-antilinear
on each simple component of $\g_1$.
Let $\lh_1$, $\lh_2$ be the $\R$-subalgebras of $\h$ such that 
$\lh_i^\perp\subset\lh_i$ and $\h_i\subset\lh_i$.
Denote as above $\l_i =\lh_i^\perp = Ker K(\cdot,\cdot)\vert_{\lh_i}$, 
$\lr_i = \g_i + \lh_i.$
Suppose that  $\phi$ is extended  to the $\R$-isomomorphism
$$\phi :\lr_1/\l_1\mapsto\lr_2/\l_2$$
such that $\phi$ preserve $K(\cdot,\cdot)$ and has no fixed points in $\h$
(condition iv)).\\
Consider the subalgebra $\W$ (see 1.2).
Then \\
1) $(\g_\R, \W, \g_\R\oplus\g_\R))$ is a Manin triple;\\
2) Every Manin triple $(\g_\R, W',\g_\R\oplus\g_\R)$ is gauge 
equivalent to $(\g_\R, \W, \g_\R\oplus\g_\R)$.\\ 
{\bf Proof}.
The complexification $\gC$ of $\g_\R$ equals to $\g\oplus\g$.
The set of simple roots of $\gC$ is $\Pi^c= \Pi\bigcup\Pi'$ where
$\Pi =\{\a_1',\ldots,\a_n'\}$ and 
$\Pi'= \{\a_1",\ldots,\a_n"\}$ are systems of symple roots of $\g$. 
The conjugation $\q$ of $\g_\R$ acts on $\Pi^c$ by $\q(\a_i)= \a_i'$
and $\q(a_i')= \a_i$.

For the subsets $\1$, $\2 \subset\Pi$ we 
consider $\Pi_i^c= \Pi_i\bigcup\Pi_i'$.
Decompose $\1$, $\2$ in the union of connected components: 
$$\1 = \bigcup_{s=1}^k \Pi_{1s},\quad
\2 = \bigcup_{s=1}^k \Pi_{2s}$$
Suppose $\phi(\Pi_{1s})= \Pi_{2s}$.
We extend $\phi$ to
$\q$-invariant map $\phi:\1^c\mapsto\2^c$,
putting $\phi (\Pi_{1s})=\Pi_{1s}$ ,if $\phi\vert_{\g_{1s}}$ is
$\C$-linear, and 
$\phi (\Pi_{1s})=\Pi_{1s}'$ ,if $\phi\vert_{\g_{1s}}$ is
$\C$-antilinear. We get $\q$-invariant map 
$\phi:\Pi_1^c\mapsto\Pi_2^c$, satisfying i), ii), iii).  
We extend $\phi$ to $\q$-invariant $\C$-linear map 
$\lr_1(\C)/\l_1(\C)$ onto $\lr_2(\C)/\l_2(\C)$,
obeying i)-v*).
By Theorem 1.2, $\W\otimes\C$ is complementary subalgebra in Manin triple
$(\gC, \W\otimes\C, \gC\oplus\gC)$. 
Considering $\q$-invariants, we see that  
$\W$ is complementary subalgebra in real Manin triple
$(\g, \W, \g\oplus\g)$. This proves 1).

On the other hand, if $(\g_\R, W, \g_\R\oplus\g_\R)$ is a real Manin triple,
then $W$ is gauge equivalent to $\W_\R$. The last subalgebra is constructed by
the two $\q$-invariant subsets  
$\Pi_i^c= \Pi_i\bigcup\Pi_i'$ and $\q$ invariant map 
$\phi:\Pi_1^c\mapsto\Pi_2^c$, obeying i)-v*).
If $\phi (\Pi_{1s})=\Pi_{1s}$ then the restriction of $\phi$ to 
$\g_{1s}(\R)= (\g_{1s}(\C))^\Gamma$ is $\C$-linear. 
In the case $\phi (\Pi_{1s})=\Pi_{1s}'$
the restriction is $\C$-antilinear.
The algebra of invariants coinsides with $Ad_r\W$.$\Box$

\end{document}